\documentclass{article}
\usepackage{graphicx} % Required for inserting images
\usepackage[utf8]{inputenc}
\usepackage{titlesec}
\titleformat{\section}
{\normalfont\scshape}
{\thesection.}{0.5em}{\centering}
\titleformat{\subsection}
{\small\scshape}
{\thesubsection.}{0.5em}{\centering}
\titleformat{\subsubsection}
{\small\scshape}
{\thesubsubsection.}{0.5em}{\centering}
\usepackage{tikz-cd}
\usepackage{tikz}
\usepackage{comment}
\usepackage{amsfonts}
\usepackage{indentfirst}
\renewcommand{\abstract}{\small{\section*{\abstractname}}}
\usepackage{amsmath}
\usepackage{amssymb}
\usepackage{amsthm}
\usepackage{mathrsfs}
\newtheorem{theorem}{Theorem}
\newtheorem*{theorem*}{Theorem}
\newtheorem{lemma}{Lemma}
\newtheorem{proposition}{Proposition}

\newtheorem{remark}{Remark}

\usepackage{stmaryrd}
\newcommand{\hooklongrightarrow}{\lhook\joinrel\longrightarrow}

\usepackage[
backend=biber,
style=alphabetic,
sorting=ynt
]{biblatex}
\title{\Large{\textbf{Relative weak compactness in infinite-dimensional Fefferman-Meyer duality}}}
\author{\normalsize{\textsc{Vasily Melnikov}}}
\date{\normalsize{\textsc{September 2024}}}

\begin{document}

\maketitle
\begin{abstract}
     Let $E$ be a Banach space such that $E'$ has the Radon-Nikodým property. The aim of this work is to connect relative weak compactness in the $E$-valued martingale Hardy space $H^{1}(\mu,E)$ to a convex compactness criterion in a weaker topology, such as the topology of uniform convergence on compacts in measure. These results represent a dynamic version of the deep result of Diestel, Ruess, and Schachermayer on relative weak compactness in $L^{1}(\mu,E)$.  In the reflexive case, we obtain a Kadec-Pełczyński dichotomy for $H^{1}(\mu,E)$-bounded sequences, which decomposes a subsequence into a relatively weakly compact part, a pointwise weakly convexly convergent part, and a null part converging to zero uniformly on compacts in measure. As a corollary, we investigate a parameterized version of the vector-valued Komlós theorem without the assumption of $H^{1}(\mu,E)$-boundedness.
\end{abstract}
\section{Introduction}\label{sec:intro}
Since the proof of Komlós's theorem \cite{ogkomlos}, the study of convex compactness (see \cite{convcomp}) and related functional-analytic techniques have been fundamental to modern probability theory, with applications to a diverse range of topics, including convex optimization (see \cite{convcomp}), the Bichteler-Dellacherie characterization of semimartingales and arbitrage theory (see \cite{arb}, or \cite{ftap}), and the Doob-Meyer decomposition (see \cite{doobmey}). Roughly speaking, although a bounded sequence in an infinite-dimensional space of random variables rarely admits a convergent subsequence (c.f. Theorem 2.1, \cite{gird}), it is possible—after passing to convex combinations—to obtain a sequence that converges almost everywhere. These results generalize to arbitrary measure spaces (see \cite{infinitekom}), and have dynamic counterparts for certain classes of stochastic processes (see \cite{martconvcomp}, or \cite{schcamp}).
\par
Similarly to this situation, given a reflexive Banach space $E$ and a bounded sequence $\{x_{n}\}_{n}\subset E$, it is a consequence of the Eberlein-Šmulian-Grothendieck theorem that it is possible to find convex combinations $y_{m}\in\mathrm{co}\{x_{n}:n\geq m\}$ such that $\{y_{m}\}_{m}$ converges in norm. More generally, a bounded subset $K\subset E$ of an arbitrary Banach space $E$ is relatively weakly compact if, and only if, the conclusion of the last sentence holds for all sequences $\{x_{n}\}_{n}\subset K$.
\par
Given how similar convex compactness results in probability theory are to the Eberlein-Šmulian-Grothendieck theorem, it is not surprising that combining the two concepts leads to a fruitful characterization of weak compactness in spaces of Bochner-integrable functions (see \cite{schruessdies}, or \cite{infinity-eberlein}) which generalizes the Dunford-Pettis criterion. In particular, a convex compactness criterion in a weaker topology (the topology of convergence in measure) and uniform integrability is both necessary and sufficient for relative weak compactness in the Bochner-Lebesgue space $L^{1}(\mu,E)$:
\begin{theorem*}{(Diestel-Ruess-Schachermayer)}
    The following are equivalent for $K\subset L^{1}(\mu,E)$:
    \begin{enumerate}
        \item $K\subset L^{1}(\mu,E)$ is relatively weakly compact.
        \item $K$ is uniformly integrable, and for each sequence $\{\xi_{n}\}_{n}\subset K$ there exists $\zeta_{n}\in\mathrm{co}\{\xi_{m}:m\geq n\}$ with $\{\zeta_{n}\}_{n}$ converging in measure.
    \end{enumerate}
\end{theorem*}
The Dellacherie-Meyer-Yor characterization of relative weak compactness in the martingale Hardy space is the analogue of the Dunford-Pettis theorem to $H^{1}$-martingales (see \cite{delmeyyor}). Thus, it is natural to inquire whether or not the above `vector-valued Dunford-Pettis theorem' generalizes to vector-valued $H^{1}$-martingales.
\par
We apply the techniques of convex compactness to characterize relative weak compactness in the $E$-valued martingale Hardy space $H^{1}(\mu,E)$, where $E$ is a Banach space such that $E'$ has the Radon-Nikodým property. In particular, it is shown that a uniform integrability condition and a sequential convex compactness requirement in the u.c.p. topology\footnote{The u.c.p. topology is the topology of uniform convergence on compacts in probability; see \S\ref{sec:prelim} for details.} are necessary and sufficient for relative weak compactness in $H^{1}(\mu,E)$ to hold. For reflexive spaces, it is shown that one can drop this convex compactness assumption. Furthermore, we investigate sequential convex compactness properties for sequences of $E$-valued martingales under the assumption that $E$ is reflexive. All results of this form are obtained using a martingale Kadec-Pełczyński dichotomy, which decomposes a bounded sequence in $H^{1}(\mu,E)$ (after passing to a subsequence) into a relatively weakly compact `regular part', and a `singular part'. The singular part either converges to zero in the u.c.p. topology, or can be further decomposed into a predictable finite variation part, and a u.c.p. null part.
\par
The structure of the paper is as follows. In Section \ref{sec:prelim}, we establish our notation. Section \ref{sec:dmy} contains the first main result of this paper, Theorem \ref{thm:del-mey-yor} on relative weak compactness in $H^{1}(\mu,E)$. Section \ref{sec:kadec} contains the second main result of the paper, Theorem \ref{thm:ucp-komlos}, which yields a Kadec-Pełczyński decomposition for $H^{1}(\mu,E)$-bounded sequences, where $E$ is reflexive. We end Section \ref{sec:kadec} by applying Theorem \ref{thm:ucp-komlos} to investigate a parameterized version of the vector-valued Komlós theorem, Theorem \ref{thm:unbound-komlos}.
\section{Preliminaries}\label{sec:prelim}
Let $E$ be a Banach space, and denote by $E'$ the dual of $E$. The norms of $E$ and $E'$ will be denoted $\Vert{\cdot}\Vert_{E}$ and $\Vert{\cdot}\Vert_{E'}$ (respectively); the closed unit balls of $E$ and $E'$ will be denoted $B_{E}$ and $B_{E'}$ (respectively). If $F$ and $G$ are vector spaces, and $\langle{\cdot,\cdot}\rangle:F\times G\longrightarrow\mathbb{R}$ is a bilinear map, then $\sigma(F,G,\langle{\cdot,\cdot}\rangle)$ denotes the locally convex topology on $F$ defined by the family of seminorms $g\longmapsto\vert{\langle{\cdot,g}\rangle}\vert$. If the duality pairing is understood (as will be the case in this article), then we will simply write $\sigma(F,G)$. In particular, $\sigma(E,E')$ coincides with the weak topology on $E$, and $\sigma(E',E)$ coincides with the weak-star topology on $E'$. Denote by $\mathscr{B}(E)$ the Borel $\sigma$-algebra generated by the norm topology, and denote by $\mathscr{B}(E')$ the Borel $\sigma$-algebra generated by the dual norm topology.
\par
If $Y$ and $Z$ are nonnegative numerical quantities (potentially depending on certain parameters), we will write $Y\lesssim Z$ to mean $Y\leq cZ$ for some constant $c\geq0$ independent of any parameters.
\par
Let $(X,\mathscr{F},\{\mathscr{F}_{t}:t\geq0\},\mu)$ be a stochastic basis, i.e., $\mathbb{F}=\{\mathscr{F}_{t}:t\geq0\}$ is a filtration of sub-$\sigma$-algebras of $\mathscr{F}$, $\mathscr{F}_{0}$ is the $\sigma$-algebra generated by the null sets $\mathscr{N}$ of $\mu$, and
\begin{equation*}
    \mathscr{F}_{t}=\bigcap_{s>t}\mathscr{F}_{s},
\end{equation*}
for all $t\geq0$. Denote by $\mathscr{P}$ the predictable $\sigma$-algebra on $[0,\infty)\times X$; we will say that an $E$-valued process $M$ is predictable if the map $(t,\omega)\longmapsto M_{t}(\omega)$ is $\mathscr{P}$-measurable when the codomain is equipped with $\mathscr{B}(E)$.
Recall that a subset $A\subset [0,\infty)\times X$ is said to be \textit{evanescent} if
\begin{equation*}
    \{\omega\in X:\exists t\in[0,\infty)\textrm{ with }(t,\omega)\in A\}\in\mathscr{N}.
\end{equation*}
Note that the countable union of evanescent sets is evanescent. In general, we will consider processes equivalent if they agree outside of an evanescent set. This contrasts with the `obvious' notion of equality between stochastic processes, equivalence up to \textit{modification}. Two $\mathbb{F}$-adapted processes $M$ and $N$ are modifications of each other if
\begin{equation*}
    \mu(\{M_{t}=N_{t}\})=1,
\end{equation*}
for each $t\in[0,\infty)$.
\par
On the space of $\mathbb{F}$-adapted $E$-valued càdlàg processes modulo evanescence, we will consider the u.c.p. (`uniform convergence on compacts in probability') topology, defined by the translation-invariant metric
\begin{equation*}
    D(M,N)=\sum_{n=1}^{\infty}\frac{1}{2^{n}}\int_{X}1\wedge\left(M-N\right)^{\ast}_{n}d\mu,
\end{equation*}
where $M^{\ast}_{t}(\cdot)=\sup_{s\leq t}\left\Vert{M_{s}(\cdot)}\right\Vert_{E}$. By (Proposition 2.2, \cite{decomp-yaro}), the space of $\mathbb{F}$-adapted $E$-valued càdlàg processes is Cauchy complete with respect to this metric. If $\{M^{n}\}_{n}\longrightarrow M$ in the u.c.p. topology, then there is a subsequence $\{n_{k}\}_{k}$ such that $\{M^{n_{k}}\}_{k}\longrightarrow M$ uniformly on all compact intervals, up to a $\mu$-null set.
\par
Recall the martingale Hardy space $H^{1}(\mu,E)$. Let $M$ be a càdlàg $E$-valued martingale. We will say that $M$ is bounded in $H^{1}(\mu,E)$, and write $M\in H^{1}(\mu,E)$, if
\begin{equation*}
    M^{\ast}_{\infty}=\sup_{t>0}M^{\ast}_{t}\in L^{1}(\mu),
\end{equation*}
and there exists a random variable $M_{\infty}\in L^{1}(\mu,E)$, the \textit{terminal value} of $M$, such that $M_{t}$ is the conditional expectation of $M_{\infty}$ with respect to $\mathscr{F}_{t}$ for each $t\in[0,\infty)$. Under the norm $M\longmapsto\Vert{M^{\ast}_{\infty}}\Vert_{L^{1}}$, $H^{1}(\mu,E)$ is a Banach space. If $E=\mathbb{R}$, we will denote $H^{1}(\mu,\mathbb{R})$ by $H^{1}(\mu)$. 
\par
In general, the condition that $M^{\ast}_{\infty}\in L^{1}(\mu)$ \textit{does not} imply that an $E$-valued martingale $M$ is the conditional expectation process of some random variable $M_{\infty}\in L^{1}(\mu,E)$. Fortunately, there are structural assumptions on $E$ which ensure the existence of a terminal value for such an $M$, namely, that $E$ possesses the Radon-Nikodým property.
\par
Given a subset $K\subset H^{1}(\mu,E)$ of $E$-valued Hardy martingales, denote by
\begin{equation*}
    K^{\ast}=\{M^{\ast}_{\infty}:M\in K\}\subset L^{1}(\mu)
\end{equation*}
the set of maximal functions, evaluated at infinity, of martingales in $K$.
\par
Recall that a subset $K\subset L^{1}(\mu)$ is said to be \textit{uniformly integrable} if $K$ is $L^{1}$-bounded and, for each $\varepsilon>0$, there exists $\delta>0$ such that
\begin{equation*}
    \sup_{\xi\in K}\int_{A}\vert{\xi}\vert d\mu\leq\varepsilon,
\end{equation*}
whenever $A\in\mathscr{F}$ satisfies $\mu(A)\leq\delta$. Vitali's convergence theorem states that a sequence $\{\xi_{n}\}_{n}\subset L^{1}(\mu)$ converges to $\xi$ in $L^{1}(\mu)$ if, and only if, $\{\xi_{n}\}_{n}$ is uniformly integrable and converges in measure to $\xi$.
\section{Relative weak compactness in $H^{1}(\mu,E)$}\label{sec:dmy}
The Dellacherie-Meyer-Yor characterization of weak compactness in $H^{1}(\mu)=H^{1}(\mu,\mathbb{R})$ is the equivalence of the following (see Théorème 1, \cite{delmeyyor}).
\begin{enumerate}
    \item $K\subset H^{1}(\mu)$ is relatively weakly compact in $H^{1}(\mu)$.
    \item $K^{\ast}=\{M^{\ast}_{\infty}:M\in K\}\subset L^{1}(\mu)$ is uniformly integrable.
    \item $\left\{[M,M]^{1/2}_{\infty}:M\in K\right\}\subset L^{1}(\mu)$ is uniformly integrable.
\end{enumerate}
\par
Suppose that one replaces $\mathbb{R}$ as above with $E$, an arbitrary Banach space. Without additional assumptions on $E$, condition (2) is not equivalent to (1). Indeed, take a non-reflexive Banach space $E$ and a bounded sequence $\{x_{n}\}_{n}\subset E$ without a weakly convergent subsequence (such a sequence exists, since $B_{E}$ is not sequentially $\sigma(E,E')$-compact); define $\{M^{n}\}_{n}\subset H^{1}(\mu,E)$ by setting $M^{n}_{t}=x_{n}$ for all $t$. Notice that $\{M^{n}\}_{n}$ cannot have a weakly convergent subsequence, despite condition (2) holding for $K=\{M^{n}:n\in\mathbb{N}\}$.
\par
Furthermore, condition (3) may not bear any relation to the other conditions. Indeed, the existence of quadratic variations combined with a suitable Burkholder-Davis-Gundy inequality is equivalent to the UMD property (see Remark 5.2, \cite{yar}).
\par
We wish to generalize the Dellacherie-Meyer-Yor theorem to vector-valued $H^{1}$-martingales. If one views the Dellacherie-Meyer-Yor criterion as an $H^{1}$-version of the Dunford-Pettis theorem, then applying the `convex compactness' philosophy of (Theorem 2.1, \cite{schruessdies}) yields the following. A set $K\subset H^{1}(\mu,E)$ should be relatively weakly compact if any sequence in $K$ is convergent in some sufficiently strong sense after passing to convex combinations, and the set of maximal functions of elements of $K$ is uniformly integrable. We will now establish a result of this form.
\begin{theorem}\label{thm:del-mey-yor}
    Suppose that $E'$ has the Radon-Nikodým property. Let $K\subset H^{1}(\mu,E)$. Then the following are equivalent.
    \begin{enumerate}
        \item $K$ is relatively weakly compact in $H^{1}(\mu,E)$.
        \item $K^{\ast}$ is uniformly integrable, and for every $\{M^{n}\}_{n}\subset K$ there exists convex combinations $N^{i}\in\mathrm{co}\{M^{n}:n\geq i\}$ and $N\in H^{1}(\mu,E)$ such that $\left\{N^{i}\right\}_{i}$ converges in the u.c.p. topology to $N$.
        \item $K^{\ast}$ is uniformly integrable, and for every $\{M^{n}\}_{n}\subset K$ there exists convex combinations $N^{i}\in\mathrm{co}\{M^{n}:n\geq i\}$ and $N\in H^{1}(\mu,E)$ such that, for each $t\in[0,\infty)$, $\left\{N^{i}_{t}\right\}_{i}$ converges $\mu$-a.e. to $N_{t}$ in $\sigma(E,E')$.
    \end{enumerate}
\end{theorem}
\begin{remark}
    The condition that $E'$ has the Radon-Nikodým property is not surprising, as this condition is used in the literature to connect relative weak compactness in $L^{1}(\mu,E)$ with a condition on the conditional expectation with respect to certain collections of sub-$\sigma$-algebras of $\mathscr{F}$ (see, for example, Theorem 1, \cite{advmath}).
\end{remark}
\begin{remark}\label{rem:rnpweaker}
    If $E$ has the Radon-Nikodým property, then we do not have to assume that $N$ from condition (2) of Theorem \ref{thm:del-mey-yor} belongs to $H^{1}(\mu,E)$. Indeed, suppose a sequence $\{M^{n}\}_{n}\subset H^{1}(\mu,E)$ is such that $\{(M^{n})^{\ast}_{\infty}\}_{n}$ is uniformly integrable, and $\{M^{n}\}_{n}$ converges in the u.c.p. topology to a process $M$. We claim that $M\in H^{1}(\mu,E)$. From uniform integrability of $\{(M^{n})^{\ast}_{\infty}\}_{n}$, Vitali's convergence theorem yields that $M$ is a martingale. Since
    \begin{equation*}
        \int_{X}M^{\ast}_{t}d\mu\leq\sup_{n}\int_{X}(M^{n})^{\ast}_{\infty}d\mu<\infty,
    \end{equation*}
    by Fatou's lemma for each $t\geq0$, $M^{\ast}_{\infty}\in L^{1}(\mu)$. All that is left to show is the existence of a terminal value for $M$, which follows from the Radon-Nikodým property of $E$.
\end{remark}
We note that the equivalence between condition (2) and (1) is certainly unexpected. Indeed, weak convergence in $H^{1}(\mu,E)$ depends also on the values of the processes at infinity, while u.c.p. convergence depends only on finite times.
\par
We will need the following lemma. Essentially, it confirms the surprising result mentioned in the previous paragraph.
\begin{lemma}\label{lem:doob-mart}
    Let $\{N^{i}\}_{i}\subset H^{1}(\mu,E)$ converge to zero in the u.c.p. topology, where $E'$ has the Radon-Nikodým property. If $\{(N^{i})^{\ast}_{\infty}\}_{i}$ is uniformly integrable, then there are convex combinations $L^{j}\in\mathrm{co}\{N^{i}:i\geq j\}$ such that $\{L^{j}\}_{j}$ converges to zero in $H^{1}(\mu,E)$.
\end{lemma}
\begin{proof}
    We may assume that $\mathscr{F}$ is generated by $\bigcup_{n=1}^{\infty}\mathscr{F}_{n}$. Fix $z'\in L^{1}(\mu,E)'$; denote by $\xi_{i}$ the terminal value of $N^{i}$. Denote $C=\sup_{i}\Vert{\xi_{i}}\Vert_{L^{1}(\mu,E)}$; if $C=0$, the claim is trivially true, so we may assume that $C>0$. By (Section IV.1, Theorem 1, \cite{vec-measures}), there exists an essentially bounded $\mathscr{B}(E')$-measurable function $y':X\longrightarrow E'$ such that
    \begin{equation*}
        \langle{\xi,z'}\rangle=\int_{X}\langle{\xi,y'}\rangle d\mu,
    \end{equation*}
    for each $\xi\in L^{1}(\mu,E)$. Our goal will be to show that $\{\vert{\langle{\xi_{i},z'}\rangle}\vert\}_{i}$ converges to zero; thus, by normalization, one may assume that $\Vert{y'\Vert}_{E'}\leq1$ up to a $\mu$-null set.
    \par
    Let $\varepsilon>0$ be arbitrary. Let $\delta>0$ be such that
    \begin{equation}\label{eq:delta1}
        \sup_{i}\int_{B}\left(N^{i}\right)^{\ast}_{\infty}d\mu<\frac{\varepsilon}{6},
    \end{equation}
    whenever $B\in\mathscr{F}$ satisfies $\mu(B)\leq\delta$. By Pettis's measurability theorem, we may assume that the range of $y'$ is separable in the dual norm topology. Thus, by Egorov's theorem, there is a simple function $u'=\sum_{i=1}^{r}x'_{i}\otimes\mathbf{1}_{A_{i}}$ with $\{x'_{1},\dots,x'_{r}\}\subset B_{E'}$ and $\{A_{i}:i\in\mathbb{N},i\leq r\}$ pairwise disjoint, such that $\Vert{\mathbf{1}_{D}(y'-u')}\Vert_{L^{\infty}(\mu,E)}<\frac{\varepsilon}{3C}$ for some measurable $D\in\mathscr{F}$ with $\mu(X\setminus D)\leq\delta$.
    \par
    Fix $\eta>0$ such that
    \begin{equation}\label{eq:eta1}
        \sup_{i}\int_{B}\left(N^{i}\right)^{\ast}_{\infty}d\mu<\frac{\varepsilon}{12r},
    \end{equation}
    whenever $\mu(B)\leq\eta$. Let $\mathscr{G}$ denote the Boolean algebra $\bigcup_{n=1}^{\infty}\mathscr{F}_{n}$; denote $\nu=\mu|_{\mathscr{G}}$. Applying the $\pi$-$\lambda$-theorem (as done in Theorem 3.3, \cite{billing}), it follows that the Carathéodory extension $\nu^{\ast}$ of $\nu$ defined by
    \begin{equation*}
        A\longmapsto\nu^{\ast}(A)=\inf\left\{\sum_{i=1}^{\infty}\nu(S_{i}):A\subset\bigcup_{n=1}^{\infty}S_{i},S_{i}\in\mathscr{G}\textrm{ for all }i\right\}
    \end{equation*}
    is such that $\nu^{\ast}|_{\mathscr{F}}=\mu$. Fix $i\in\mathbb{N}$, $i\leq r$. Thus, there exists a pairwise disjoint sequence $\{S_{n}^{i}\}_{n}\subset\mathscr{G}$ such that $A_{i}\subset\bigcup_{n=1}^{\infty}S_{n}^{i}$, and
    \begin{equation*}
        \sum_{n=1}^{\infty}\mu(S^{i}_{n})<\mu(A_{i})+\eta,
    \end{equation*}
    implying $\mu\left((X\setminus A_{i})\cap\bigcup_{n=1}^{\infty}S_{n}^{i}\right)<\eta$. Since $\sum_{n=1}^{\infty}\mu(S^{i}_{n})<\infty$ there exists $\ell_{i}\in\mathbb{N}$, $\ell_{i}\geq 2$ such that
    \begin{equation*}
        \sum_{n=\ell_{i}}^{\infty}\mu(S^{i}_{n})<\eta.
    \end{equation*}
    Let $k\in\mathbb{N}$ be such that $\{S_{1}^{i},\dots,S_{\ell_{i}-1}^{i}\}\subset\mathscr{F}_{k}$ for each $i\in\mathbb{N}$, $i\leq r$. Such a $k$ exists; indeed, $r<\infty$. Let $j\in\mathbb{N}$ be such that
    \begin{equation}\label{eq:loc-h1-convg}
        \int_{X}\left(N^{i}\right)^{\ast}_{k}d\mu<\frac{\varepsilon}{6r},
    \end{equation}
    whenever $i\geq j$; such a $j$ exists by Vitali's convergence theorem (indeed, $\{N^{i}\}_{i}$ converges to zero in the u.c.p. topology).
    \par
    For each $i$,
    \begin{equation*}
        \vert{\langle{\xi_{i},z'}\rangle}\vert=\left\vert{\int_{X}\langle{\xi_{i},y'}\rangle}d\mu\right\vert\leq\left\vert{\int_{X}\langle{\xi_{i},y'-u'}\rangle}d\mu\right\vert+\left\vert{\int_{X}\langle{\xi_{i},u'}\rangle}d\mu\right\vert
    \end{equation*}
    \begin{equation*}
        \leq\left\vert{\int_{D}\langle{\xi_{i},y'-u'}\rangle}d\mu\right\vert+\left\vert{\int_{X\setminus D}\langle{\xi_{i},y'}\rangle}d\mu\right\vert+\left\vert{\int_{X\setminus D}\langle{\xi_{i},u'}\rangle}d\mu\right\vert+\left\vert{\int_{X}\langle{\xi_{i},u'}\rangle}d\mu\right\vert
    \end{equation*}
    \begin{equation}\label{eq:global-bound}
        <\frac{\varepsilon}{3C}\int_{X}\left\Vert{\xi_{i}}\right\Vert_{E}d\mu+\frac{\varepsilon}{6}+\frac{\varepsilon}{6}+\left\vert{\int_{X}\langle{\xi_{i},u'}\rangle}d\mu\right\vert\leq\frac{2\varepsilon}{3}+\left\vert{\int_{X}\langle{\xi_{i},u'}\rangle}d\mu\right\vert,
    \end{equation}
    by the triangle inequality and (\ref{eq:delta1}), noting that $\mu(X\setminus D)\leq\delta$. Furthermore,
    \begin{equation*}
        \left\vert{\int_{X}\langle{\xi_{i},u'}\rangle}d\mu\right\vert\leq\sum_{q=1}^{r}\left\vert\left\langle{\int_{A_{q}}\xi_{i}d\mu,x'_{q}}\right\rangle\right\vert
    \end{equation*}
    \begin{equation*}
        =\sum_{q=1}^{r}\left\vert\left\langle{\int_{\bigcup_{n=1}^{\infty}S^{q}_{n}}\xi_{i}d\mu-\int_{\left(\bigcup_{n=1}^{\infty}S^{q}_{n}\right)\cap(X\setminus A_{q})}\xi_{i}d\mu,x'_{q}}\right\rangle\right\vert
    \end{equation*}
    \begin{equation*}
        \leq\sum_{q=1}^{r}\left(\left\vert\left\langle{\int_{\bigcup_{n=1}^{\infty}S^{q}_{n}}\xi_{i}d\mu,x'_{q}}\right\rangle\right\vert+\left\vert\left\langle{\int_{\left(\bigcup_{n=1}^{\infty}S^{q}_{n}\right)\cap(X\setminus A_{q})}\xi_{i}d\mu,x'_{q}}\right\rangle\right\vert\right)
    \end{equation*}
    \begin{equation}\label{eq:local-bound}
        \leq\sum_{q=1}^{r}\left\vert\left\langle{\int_{\bigcup_{n=1}^{\infty}S^{q}_{n}}\xi_{i}d\mu,x'_{q}}\right\rangle\right\vert+\frac{\varepsilon}{12},
    \end{equation}
    by the triangle inequality and (\ref{eq:eta1}). If $i\geq j$,
    \begin{equation*}
        \left\vert\left\langle{\int_{\bigcup_{n=1}^{\infty}S^{q}_{n}}\xi_{i}d\mu,x'_{q}}\right\rangle\right\vert\leq\left\vert\left\langle{\int_{\bigcup_{n=1}^{\ell_{i}-1}S^{q}_{n}}\xi_{i}d\mu,x'_{q}}\right\rangle\right\vert+\left\vert\left\langle{\int_{\bigcup_{n=\ell_{i}}^{\infty}S^{q}_{n}}\xi_{i}d\mu,x'_{q}}\right\rangle\right\vert
    \end{equation*}
    \begin{equation*}
        \leq\left\vert\left\langle{\int_{\bigcup_{n=1}^{\ell_{i}-1}S^{q}_{n}}N^{i}_{k}d\mu,x'_{q}}\right\rangle\right\vert+\frac{\varepsilon}{12r}<\frac{\varepsilon}{4r}
        \end{equation*}
    by (\ref{eq:eta1}) (noting that $\mu(\bigcup_{n=\ell_{i}}^{\infty}S^{q}_{n})<\eta$), (\ref{eq:loc-h1-convg}), and the martingale property, as $\bigcup_{n=1}^{\ell_{i}-1}S^{q}_{n}\in\mathscr{F}_{k}$ for all $q\in\mathbb{N}$ with $q\leq r$. Combining with (\ref{eq:local-bound}), one obtains
    \begin{equation*}
        \left\vert{\int_{X}\langle{\xi_{i},u'}\rangle}d\mu\right\vert<\frac{\varepsilon}{12}+\sum_{q=1}^{r}\frac{\varepsilon}{4r}=\frac{\varepsilon}{3},
    \end{equation*}
    if $i\geq j$. Combining with (\ref{eq:global-bound}), this shows that
    \begin{equation*}
        \vert{\langle{\xi_{i},z'}\rangle}\vert<\varepsilon,
    \end{equation*}
    if $i\geq j$. Since $z'$ and $\varepsilon$ were arbitrary, this shows that $\{\xi_{i}\}_{i}$ converges to zero in the weak topology on $L^{1}(\mu,E)$.
    \par
    Thus, by Mazur's lemma, there exists convex combinations $L^{j}\in\mathrm{co}\{N^{i}:i\geq j\}$ such that the terminal values of $\{L^{j}\}_{j}$ converge to zero in the norm topology on $L^{1}(\mu,E)$. By Doob's maximal inequality, $\{(L^{i})^{\ast}_{\infty}\}_{i}$ converges to zero in measure; by Vitali's convergence theorem, $\{(L^{i})^{\ast}_{\infty}\}_{i}$ therefore converges to zero in $L^{1}(\mu)$. Thus
    \begin{equation*}
        \lim_{j\to\infty}\Vert{L^{j}}\Vert_{H^{1}(\mu,E)}=0,
    \end{equation*}
    as desired.
\end{proof}
\begin{remark}\label{rem:tala}
    The argument in Lemma \ref{lem:doob-mart} relies heavily on $E'$ possessing the Radon-Nikodým property. Indeed, if the natural embedding $L^{\infty}(\mu,E')\hooklongrightarrow L^{1}(\mu,E)'$ is an isomorphism, then $E'$ must have the Radon-Nikodým property with respect to $\mu$ (see Section IV.1, Theorem 1, \cite{vec-measures}).
\end{remark}
\begin{proof}[Proof of Theorem \ref{thm:del-mey-yor}]
    We will first show that (1) implies (2). By the Eberlein-Šmulian-Grothendieck theorem (see Corollary 2.2, \cite{schruessdies}), for every sequence $\left\{M^{n}\right\}_{n}\subset K$, there exists $N^{i}\in\mathrm{co}\{M^{n}:n\geq i\}$ such that $\{N^{i}\}_{i}$ converges in $H^{1}(\mu,E)$ to some $N\in H^{1}(\mu,E)$. In particular, $\{(N-N^{i})^{\ast}_{\infty}\}_{i}$ tends to zero in $L^{1}(\mu)$. Applying Markov's inequality shows that $\{N^{i}\}_{i}$ converges in the u.c.p. topology to $N\in H^{1}(\mu,E)$. It remains to show that $K^{\ast}$ is uniformly integrable. Note that $K^{\ast}$ is bounded in $L^{1}(\mu)$ by the uniform boundedness principle, so (Lemme 5, \cite{delmeyyor}) implies it suffices to show that
    \begin{equation*}
        \left\{\Vert{M_{T}}\Vert_{E}:M\in K\right\}
    \end{equation*}
    is uniformly integrable for each nonnegative random variable $T$. The assignment $M\longmapsto M_{T}$ defines a continuous linear operator from $H^{1}(\mu,E)$ to $L^{1}(\mu,E)$, and therefore maps weakly compact sets to weakly compact sets. By weak compactness of $K$ in $H^{1}(\mu,E)$, $\{M_{T}:M\in K\}$ is therefore weakly compact in $L^{1}(\mu,E)$, and so $\left\{\Vert{M_{T}}\Vert_{E}:M\in K\right\}$ is uniformly integrable by (Theorem 2.1, \cite{schruessdies}). Thus, (1) implies (2).
\par
We will now show that (2) implies (1). Let $\{M^{n}\}_{n}\subset K$; we may pass to forward convex combinations $\{N^{i}\}_{i}$ such that $\{N^{i}\}_{i}$ converges to an $H^{1}(\mu,E)$-martingale $N$ in the u.c.p. topology. Thus, $\{N^{i}-N\}_{i}$ converges to zero in the u.c.p. topology. Applying Lemma \ref{lem:doob-mart} yields convex combinations $L^{j}\in\mathrm{co}\{N^{i}-N:i\geq j\}$ such that $\{L^{j}\}_{j}$ converges to zero in $H^{1}(\mu,E)$. Clearly, $L^{j}+N\in\mathrm{co}\{N^{i}:i\geq j\}$ and converges in $H^{1}(\mu,E)$ to $N$. Applying the Eberlein-Šmulian-Grothendieck theorem shows that $K$ is relatively weakly compact in $H^{1}(\mu,E)$.
\par
It remains to show that (3) is equivalent to (2). Since u.c.p. convergence—by virtue of the Borel-Cantelli lemma and diagonalization—implies the existence of a subsequence which converges uniformly on compacts $\mu$-a.e., (2) implies (3). For the reverse implication, let $\{M^{n}\}_{n}\subset K$ be arbitrary. Pass to convex combinations $N^{i}\in\mathrm{co}\{M^{n}:n\geq i\}$ converging in the sense of condition (3) to $N\in H^{1}(\mu,E)$; applying the same argument as (Theorem 2.1, \cite{schruessdies}) and passing to convex combinations, one may assume that $\{N^{i}_{t}\}_{i}$ converges in $L^{1}(\mu,E)$ to $N_{t}$ for each $t\geq0$. By Doob's maximal inequality
\begin{equation*}
    \mu\left(\left\{(N-N^{i})^{\ast}_{t}\geq\varepsilon\right\}\right)\leq\frac{1}{\varepsilon}\int_{X}\left\Vert{N_{t}-N^{i}_{t}}\right\Vert_{E}d\mu,
\end{equation*}
and so $\{N^{i}\}_{i}$ converges to $N\in H^{1}(\mu,E)$ in the u.c.p. topology, as desired.
\end{proof}
\par
Let us note that for non-reflexive Banach spaces $E$, conditions (2) or (3) of Theorem \ref{thm:del-mey-yor} will not be satisfied for all subsets $K\subset H^{1}(\mu,E)$ such that $K^{\ast}$ is uniformly integrable. However, in the reflexive case, the compactness criterion given by Theorem \ref{thm:del-mey-yor} is exactly the analogous version of the Dellacherie-Meyer-Yor criterion. Indeed, one has the following result.
\begin{proposition}\label{prop:reflex-del-mey-yor}
    Suppose that $E$ is reflexive. The following are equivalent for $K\subset H^{1}(\mu,E)$.
    \begin{enumerate}
        \item $K$ is relatively weakly compact in $H^{1}(\mu,E)$.
        \item $K^{\ast}$ is uniformly integrable.
    \end{enumerate}
\end{proposition}
\begin{proof}
    The same technique used in the proof of Theorem \ref{thm:del-mey-yor} shows that (1) implies (2). We will now show the converse.
    \par
    By the Eberlein-Šmulian-Grothendieck theorem, it suffices to show that for every $\{M^{n}\}_{n}\subset K$ there exists convex combinations $N^{i}\in\mathrm{co}\{M^{n}:n\geq i\}$ such that $\{N^{i}\}_{i}$ converges in $H^{1}(\mu,E)$. Using the `vulgar' version of the vector-valued Komlós theorem (see Theorem 1.4, \cite{del-sch}), we may find convex combinations $N^{i}\in\mathrm{co}\{M^{n}:n\geq i\}$ such that the corresponding terminal values $\{N^{i}_{\infty}\}_{i}$ converge $\mu$-a.e. in the norm topology to some $\zeta\in L^{1}(\mu,E)$. Using a similar argument to Remark \ref{rem:rnpweaker}, we may find $N\in H^{1}(\mu,E)$ with terminal value $\zeta$. Uniform integrability of $K^{\ast}$ combined with Vitali's convergence theorem shows that $\{N^{i}_{\infty}\}_{i}$ is convergent in $L^{1}(\mu,E)$ to $N_{\infty}$. By Doob's maximal inequality, $\{(N-N^{i})^{\ast}_{\infty}\}_{i}$ converges to zero in probability, hence in $L^{1}(\mu)$ by Vitali's convergence theorem and uniform integrability of $K^{\ast}$. Thus, $\{N^{i}\}_{i}$ converges in $H^{1}(\mu,E)$, as desired.
\end{proof}
We will now address to what extent Theorem \ref{thm:del-mey-yor} can be generalized to arbitrary Banach spaces. As Remark \ref{rem:tala} shows, the proof of Lemma \ref{lem:doob-mart} depends on $E'$ possessing the Radon-Nikodým property with respect to $\mu$. However, the exact argument of Lemma \ref{lem:doob-mart} still shows that convergence to zero holds in the weak topology on the dual pair $\langle{L^{1}(\mu,E),L^{\infty}(\mu,E')}\rangle$. We refer the reader to \cite{studiacomp} for sufficient conditions for this to imply weak or norm convergence in $L^{1}(\mu,E)$.
\par
To complement the above paragraph, we now give an example to show why Theorem \ref{thm:del-mey-yor} may fail to hold if $E'$ does not possess the Radon-Nikodým property. We work in discrete time, trusting that the reader understands how to extend to continuous time. Consider $(X,\mathscr{F},\mu)$ as the Cantor group $\{-1,1\}^{\mathbb{N}}$ equipped with its natural Haar probability measure $\mu$, with $\mathscr{F}$ the $\mu$-completion of the Borel $\sigma$-algebra of the product topology. Let the filtration $\{\mathscr{F}_{n}:n\in\mathbb{N}\cup\{0\}\}$ be defined so that $\mathscr{F}_{0}$ the $\mu$-completion of the trivial $\sigma$-algebra, and (for $n\geq 1$) $\mathscr{F}_{n}$ is the $\mu$-completion of the $\sigma$-algebra generated by $\{\varepsilon_{1},\dots,\varepsilon_{n}\}$, where $\varepsilon_{m}$ denotes the projection of $X=\{-1,1\}^{\mathbb{N}}$ onto the $m$-th coordinate (a Rademacher random variable). Let $E$ be a Banach space containing a sequence $\{x_{n}\}_{n}\subset B_{E}$ equivalent to the $\ell^{1}$ basis.\footnote{A sequence $\{x_{n}\}_{n}\subset E$ is said to be equivalent to the $\ell^{1}$ basis if,
\begin{equation*}
    \sum_{i=1}^{\infty}\vert{a_{i}}\vert\lesssim\left\Vert{\sum_{i=1}^{\infty}a_{i}x_{i}}\right\Vert_{E},
\end{equation*}
whenever $\{a_{i}\}_{i}\subset\mathbb{R}$ is a sequence with finitely many non-zero terms; essentially, if $\{x_{n}\}_{n}\subset E$ is equivalent to the $\ell^{1}$ basis, then it is as far from being relatively weakly compact as possible.} For each $n\in\mathbb{N}$, define a martingale $M^{n}\in H^{1}(\mu,E)$ by
\begin{equation*}
    M^{n}_{\infty}=\varepsilon_{n}x_{n}.
\end{equation*}
$\{M^{n}\}_{n}$ forms a bounded sequence in $H^{1}(\mu,E)$. For each $n<m$, $M^{m}_{n}=0$ $\mu$-a.e., so the maximal functions $\{(M^{n})^{\ast}_{m}\}_{n}$ converge to zero $\mu$-a.e. for each $m\in\mathbb{N}$. However, in contrast to the conclusion of Theorem \ref{thm:del-mey-yor}, the sequence $\{M^{n}\}_{n}\subset H^{1}(\mu,E)$ has no weakly convergent subsequence. Indeed, the sequence of terminal values $\{M^{n}_{\infty}\}_{n}$ is an $\ell^{1}$ basis in $L^{1}(\mu,E)$, and therefore does not admit a weakly convergent subsequence.
\section{Martingale Kadec-Pełczyński dichotomy}\label{sec:kadec}
Let $\{\xi_{n}\}_{n}\subset L^{1}(\mu)$ be bounded. By the Kadec-Pełczyński decomposition theorem \cite{ogkp}, it is possible to split (after passing to a subsequence) $\{\xi_{n}\}_{n}$ into the sum of a uniformly integrable sequence (the `regular part'), and a sequence which converges to zero $\mu$-a.e. (the `singular part').
\par
The theorem of Komlós \cite{ogkomlos} is a fundamental result in probability; functional analysts will quickly notice that the Kadec-Pełczyński decomposition theorem is a direct improvement on Komlós's theorem. Importantly, unlike Komlós's theorem, the Kadec-Pełczyński decomposition immediately provides deep and fundamental insight into the functional-analytical structure of $L^{1}(\mu)$-bounded sequences. Indeed, the Kadec-Pełczyński decomposition allows one to conclude that the $\ell^{1}$ basis is universal among non-uniformly integrable bounded sequences; after throwing away a relatively weakly compact `regular part', one is left with something that is equivalent to the $\ell^{1}$ basis (due to Rosenthal's $\ell^{1}$ theorem), and converges to the same pointwise almost everywhere limit.
\par
$H^{1}(\mu)$ versions of the Kadec-Pełczyński decomposition have been investigated (see \cite{del-sch}). However, without the assumption of continuity or additional regularity, the singular part of this decomposition may not tend to zero in any reasonable sense. This phenomenon is caused by the underlying martingales having large jumps, leading to the singular parts converging, in the sense of Fatou (see Definition 5.2, \cite{optdecompconstrain}),\footnote{Unfortunately, the notion of a Fatou limit does not generalize to the vector-valued context, so we will be forced in this section to consider a different notion of convergence.} to a non-zero finite-variation process.
\par
In this section, we develop a Kadec-Pełczyński type decomposition for $H^{1}(\mu,E)$-bounded sequences, where $E$ is a reflexive Banach space. As expected from the $E=\mathbb{R}$ case, the singular parts of the decomposition do not converge to zero in any reasonable sense, unless additional assumptions are fulfilled (see Proposition \ref{prop:banach-saks-cont}). We decompose the singular part into two further sequences: a u.c.p. null sequence, and a sequence of predictable processes whose variations are bounded in $L^{1}(\mu)$. As a corollary, we establish a version of Komlós's theorem for martingales without the assumption of $H^{1}$-boundedness.
\par
We will assume, unless otherwise stated, that $E$ stands for a reflexive Banach space.
\begin{proposition}\label{prop:banach-saks-cont}
    Let $\{M^{n}\}_{n}$ be a bounded sequence in $H^{1}(\mu,E)$, where $E$ is a reflexive Banach space. Suppose that at least one of the following conditions holds.
    \begin{enumerate}
        \item The set $\{M^{n}_{T}:n\in\mathbb{N},T\textrm{ a stopping time}\}$ is uniformly integrable.
        \item The set $\{\Delta M^{n}_{T}:n\in\mathbb{N},T\textrm{ a stopping time}\}$ is uniformly integrable.
    \end{enumerate}
    Then, there is a subsequence $\{n_{k}\}_{k}$ such that
    \begin{equation*}
        M^{n_{k}}=N^{k}+L^{k},
    \end{equation*}
    for each $k\in\mathbb{N}$, where $\{N^{k}\}_{k}\subset H^{1}(\mu,E)$ is a relatively weakly compact sequence in $H^{1}(\mu,E)$, and $\{L^{k}\}_{k}$ is a null sequence in the u.c.p. topology.
\end{proposition}
\begin{proof}
     The proof is essentially a vector-valued version of (Theorem A, \cite{del-sch}). 
     \par
     We will use a version of the Kadec-Pełczyński decomposition theorem (see Theorem 2.1, \cite{del-sch}). From the Kadec-Pełczyński theorem, there exists a strictly positive sequence $\{a_{n}\}_{n}$ converging to $\infty$ in the Alexandrov compactification $\alpha[0,\infty)=[0,\infty)\cup\{\infty\}$ of $[0,\infty)$ such that, after passing to a subsequence if necessary, $\{\left(M^{n}\right)^{\ast}_{\infty}\wedge a_{n}\}_{n}$ is uniformly integrable; by passing to a further subsequence, one may assume that
     \begin{equation*}
         \sum_{n=1}^{\infty}\frac{1}{a_{n}}<\infty.
     \end{equation*}
     Define
     \begin{equation*}
         S_{n}=\inf\left\{t:\left\Vert{M^{n}_{t}}\right\Vert_{E}\geq a_{n}\right\}.
     \end{equation*}
     Let $T_{n}=\bigwedge_{k=n}^{\infty}S_{k}$. It is easy to see from Markov's concentration inequality that $\mu(\{S_{n}<\infty\})\lesssim\frac{1}{a_{n}}$, and $\mu(\{T_{n}<\infty\})\lesssim\sum_{i=n}^{\infty}\frac{1}{a_{n}}$.
     \par
     Suppose that condition (1) holds. Then
     \begin{equation*}
         \left((M^{n})^{T_{n}}\right)^{\ast}_{\infty}\leq\left((M^{n})^{\ast}_{\infty}\wedge a_{n}\right)\vee\Vert{M^{n}_{T_{n}}}\Vert_{E},
     \end{equation*}
     showing that $\{(M^{n})^{T_{n}}\}_{n}$ is relatively weakly compact in $H^{1}(\mu,E)$ by Proposition \ref{prop:reflex-del-mey-yor}. The remaining part goes to zero in the u.c.p. topology. Likewise, if condition (2) is satisfied,
     \begin{equation*}
         \left((M^{n})^{T_{n}}\right)^{\ast}_{\infty}\leq\left((M^{n})^{\ast}_{\infty}\wedge a_{n}\right)+\Vert{\Delta M^{n}_{T_{n}}}\Vert_{E},
     \end{equation*}
     showing that $\{(M^{n})^{T_{n}}\}_{n}$ is relatively weakly compact in $H^{1}(\mu,E)$ by Proposition \ref{prop:reflex-del-mey-yor}. The remaining part goes to zero in the u.c.p. topology.
\end{proof}
As a result of Proposition \ref{prop:banach-saks-cont}, it is possible to obtain a `parameterized' Komlós theorem, as long as any of the conditions in Proposition \ref{prop:banach-saks-cont} are satisfied. Of course, it is not clear what the proper analogue of Komlós's theorem is in this context. Indeed, there are several modes of convergence that could replace $\mu$-a.e. convergence. For example, one can consider u.c.p. convergence, or pointwise convergence on an evanescent subset of $[0,\infty)\times X$. The former is addressed in Proposition \ref{prop:kom-ucp} under certain regularity assumptions (c.f. Theorem \ref{thm:unbound-komlos}), while the latter is the content of Theorem \ref{thm:ucp-komlos}.
\begin{proposition}\label{prop:kom-ucp}
    Let $\{M^{n}\}_{n}$ be a bounded sequence in $H^{1}(\mu,E)$, where $E$ is a reflexive Banach space. Suppose that at least one of the following conditions holds.
    \begin{enumerate}
        \item The set $\{M^{n}_{T}:n\in\mathbb{N},T\textrm{ a stopping time}\}$ is uniformly integrable.
        \item The set $\{\Delta M^{n}_{T}:n\in\mathbb{N},T\textrm{ a stopping time}\}$ is uniformly integrable.
    \end{enumerate}
    Then there exists $N^{k}\in\mathrm{co}\{M^{n}:n\geq k\}$ such that $\{N^{k}\}_{k}$ converges in the u.c.p. topology.
\end{proposition}
\begin{remark}
    It is possible to replace the convex combinations in Proposition \ref{prop:kom-ucp} with something more tractable. For example, suppose that $L^{2}(\mu,E)$ has the Banach-Saks property, and condition (1) or (2) of Proposition \ref{prop:kom-ucp} holds. Then, using a theorem of \cite{bourgainkomlos}, it is possible to extract a subsequence that converges in the u.c.p. topology in Cesàro mean. This result is closer to the original formulation of Komlós's theorem by \cite{ogkomlos}.
\end{remark}
It is unclear whether Proposition \ref{prop:banach-saks-cont} (and hence also Proposition \ref{prop:kom-ucp}) extends to general bounded sequences in $H^{1}(\mu,E)$. The chief issue is that the jumps $\Delta M^{n}_{T_{n}}$ may be large. The idea (due to \cite{del-sch}) to remedy this situation is to `throw away' all the jumps at time $T_{n}$, and then compensate the resulting process (so that the resulting process is still a martingale). Of course, one must pay a price: in general, this finite variation part need not converge to zero in \textit{any} sense. Furthermore, one does not in general have u.c.p. convergence.
\par
If $N$ is an $E$-valued process of finite variation, denote by $\mathrm{var}(N)$ the variation process of $N$.
\begin{theorem}\label{thm:ucp-komlos}
    Let $\{M^{n}\}_{n}$ be a bounded sequence in $H^{1}(\mu,E)$, where $E$ is a reflexive Banach space. There is a subsequence $\{n_{k}\}_{k}$ such that
    \begin{equation*}
        M^{n_{k}}=N^{k}+L^{k}+R^{k},
    \end{equation*}
    for each $k\in\mathbb{N}$, where $\{N^{k}\}_{k}\subset H^{1}(\mu,E)$ is a relatively weakly compact sequence in $H^{1}(\mu,E)$, $\{L^{k}\}_{k}$ is a null sequence in the u.c.p. topology, $\{R^{k}\}_{k}$ is predictable and satisfies
    \begin{equation*}
        \sup_{k}\int_{X}\mathrm{var}(R^{k})_{\infty}d\mu<\infty,
    \end{equation*}
    and there exists a predictable process $\widetilde{R}$ of integrable variation such that, for some $X_{0}\in\mathscr{F}$ of full measure,
    \begin{equation*}
        \lim_{k\to\infty}\widetilde{R}_{t}^{k}=\widetilde{R}_{t},
    \end{equation*}
    on $X_{0}$ for all $t\geq0$, for some sequence $\widetilde{R}^{k}\in\mathrm{co}\{R^{n}:n\geq k\}$, and where the limit is in $\sigma(E,E')$.
\end{theorem}
We will assume until the end of the proof of Theorem \ref{thm:ucp-komlos} that $E$ stands for a reflexive Banach space. Of course, this implies that both $E$ and $E'$ have the Radon-Nikodým property.
\par
Fix $n\in\mathbb{N}$, and $s\in\mathbb{Q}_{+}$. By Pettis's measurability theorem, there is a separable subspace $E_{n,s}\subset E$ such that $M^{n}_{s}\in E_{n,s}$ up to a $\mu$-null set. By right continuity,
\begin{equation*}
    M^{n}_{t}(\omega)\in\overline{\bigcup_{s\in\mathbb{Q}_{+}}\bigcup_{m=1}^{\infty}E_{m,s}},
\end{equation*}
for all $n\in\mathbb{N}$ and $t\in[0,\infty)$, for all $\omega\in X_{0}$ for some subset $X_{0}\in\mathscr{F}$ of full $\mu$-measure. Thus, if we are only dealing with sequences (or convex combinations thereof), one may assume that $E$ is separable. Indeed, all of the properties mentioned in the last paragraph pass to closed linear subspaces.
\par
We will need the following lemma, which holds in any Banach space $E$ with the Radon-Nikodým property.
\begin{lemma}\label{lem:varbound}
    Suppose $E$ is a Banach space with the Radon-Nikodým property. Let $T$ be a stopping time, and let $U\in L^{1}(\mu,E)$ be $\mathscr{F}_{T}$-measurable. Then the process $\mathbf{1}_{\llbracket{T,\infty}\llbracket}U$ decomposes as $\mathbf{1}_{\llbracket{T,\infty}\llbracket}U=M+A$, where $M$ is an $E$-valued càdlàg martingale, and $A$ is a predictable process of finite variation. Furthermore, one has the bound
    \begin{equation*}
        \int_{X}\mathrm{var}(A)_{\infty}\leq\Vert{U}\Vert_{L^{1}(\mu,E)}.
    \end{equation*}
\end{lemma}
\begin{proof}
    We refer readers to the proof of (Theorem 2.20, \cite{handbookDoobMeyer}).
\end{proof}
The following lemma holds for any Banach space $E$ which is reflexive. It will allow us to obtain a convex compactness result for the singular parts of our decomposition.
\begin{lemma}\label{lem:pred-comp}
    Suppose $E$ is a reflexive Banach space. Let $\{M^{n}\}_{n}$ be a sequence of $E$-valued predictable processes of finite variation such that
    \begin{equation*}
        \sup_{n}\int_{X}\mathrm{var}(M^{n})_{\infty}d\mu<\infty.
    \end{equation*}
    There exists $\widetilde{M}^{k}\in\mathrm{co}\{M^{n}:n\geq k\}$ such that there exists a predictable process $\widetilde{M}$ of integrable variation such that, for some $X_{0}\in\mathscr{F}$ of full measure,
    \begin{equation*}
        \lim_{k\to\infty}\widetilde{M}_{t}^{k}=\widetilde{M}_{t},
    \end{equation*}
    on $X_{0}$ for all $t\geq0$, where the limit is in $\sigma(E,E')$.
\end{lemma}
\begin{proof}
    By using (Lemma, 2.5, \cite{schruessdies}) and passing to convex combinations, one may assume that
    \begin{equation}\label{eq:sup-a.e.-bounded}
        \sup_{n}\sup_{t\geq0}\left\Vert{M^{n}_{t}}\right\Vert_{E}\leq\sup_{n}\mathrm{var}(M^{n})_{\infty}<\infty,
    \end{equation}
    up to a $\mu$-null set. Let $D\subset E'$ be a countable, norm dense subset of $E'$; such a subset exists by separability of $E$, and the assumption of reflexivity. Let $F=\mathrm{span}(D)$. By the parameterized Helly selection theorem (see Proposition 13, \cite{schcamp}) and diagonalization, there exists $\widetilde{M}^{k}\in\mathrm{co}\{M^{n}:n\geq k\}$ and $X_{0}\in\mathscr{F}$ such that $\{\widetilde{M}_{t}^{k}\}_{k}$ is Cauchy on $X_{0}$ for the uniformity defined by the family of seminorms $x'\longmapsto\vert{\langle{\cdot,x'}\rangle}\vert$ as $x'$ ranges over $F$ for each $t\geq0$. By (\ref{eq:sup-a.e.-bounded}), and the coincidence of $\sigma(E,F)$ and $\sigma(E,E')$ on bounded sets, it follows that $\{\widetilde{M}_{t}^{k}\}_{k}$ is weakly Cauchy on $X_{0}$ for all $t\geq0$. By weak sequential completness, one may assume that $\{\widetilde{M}^{k}\}_{k}$ converges pointwise on $[0,\infty)\times X_{0}$ to a predictable process $\widetilde{M}$.
    \par
    We will now show that $\widetilde{M}$ has integrable variation. For each $n\in\mathbb{N}$, denote $F_{n}=\bigoplus_{i=1}^{n}E$, $F'_{n}=\bigoplus_{i=1}^{n}E'$. By the definition of $\mathrm{var}\left(\widetilde{M}\right)$ and lower semicontinuity of the function $F_{n}\ni(x_{1},\dots,x_{n})\longmapsto\sum_{i=1}^{n}\Vert{x_{i}}\Vert_{E}$ in $\sigma(F_{n},F'_{n})$ for each $n\in\mathbb{N}$,
    \begin{equation*}
        \mathrm{var}\left(\widetilde{M}\right)_{\infty}=\sup_{n}\sup_{t_{0}<\dots<t_{n}}\sum_{i=1}^{n}\Vert{\widetilde{M}_{t_{i}}-\widetilde{M}_{t_{i-1}}}\Vert_{E}\leq\sup_{n}\sup_{t_{0}<\dots<t_{n}}\liminf_{j\to\infty}\sum_{i=1}^{n}\Vert{\widetilde{M}_{t_{i}}^{j}-\widetilde{M}_{t_{i-1}}^{j}}\Vert_{E}
    \end{equation*}
    \begin{equation*}
        \leq\liminf_{j\to\infty}\sup_{n}\sup_{t_{0}<\dots<t_{n}}\sum_{i=1}^{n}\Vert{\widetilde{M}_{t_{i}}^{j}-\widetilde{M}_{t_{i-1}}^{j}}\Vert_{E}=\liminf_{j\to\infty}\mathrm{var}\left(\widetilde{M}^{j}\right)_{\infty},
    \end{equation*}
    up to a $\mu$-null set, so that
    \begin{equation*}
        \int_{X}\mathrm{var}\left(\widetilde{M}\right)_{\infty}d\mu\leq\int_{X}\liminf_{n\to\infty}\mathrm{var}\left(\widetilde{M}^{n}\right)_{\infty}d\mu\leq\liminf_{n\to\infty}\int_{X}\mathrm{var}\left(\widetilde{M}^{n}\right)_{\infty}d\mu<\infty,
    \end{equation*}
    by Fatou's lemma, which shows that $\widetilde{M}$ has integrable variation.
\end{proof}
\begin{remark}
    Reflexivity plays a crucial role in the proof of Lemma \ref{lem:pred-comp}. Indeed, if $E$ is such that each bounded sequence $\{x_{n}\}_{n}\subset E$, interpreted as a constant stochastic process, admits convex combinations converging weakly, then the Eberlein-Šmulian-Grothendieck theorem implies the unit ball of $E$ must be weakly compact. Thus, $E$ must be reflexive.
\end{remark}
\begin{proof}[Proof of Theorem \ref{thm:ucp-komlos}.]
    From the Kadec-Pełczyński theorem, there exists a strictly positive sequence $\{a_{n}\}_{n}$ converging to $\infty$ in $\alpha[0,\infty)$ such that, after passing to a subsequence if necessary, $\{\left(M^{n}\right)^{\ast}_{\infty}\wedge a_{n}\}_{n}$ is uniformly integrable; by passing to a further subsequence, one may assume that
     \begin{equation*}
         \sum_{n=1}^{\infty}\frac{1}{a_{n}}<\infty.
     \end{equation*}
     Define
     \begin{equation*}
         S_{n}=\inf\left\{t:\left\Vert{M^{n}_{t}}\right\Vert_{E}\geq a_{n}\right\}.
     \end{equation*}
     Let $T_{n}=\bigwedge_{k=n}^{\infty}S_{k}$; it is easy to see from Markov's concentration inequality that $\mu(\{S_{n}<\infty\})\lesssim\frac{1}{a_{n}}$, and $\mu(\{T_{n}<\infty\})\lesssim\sum_{i=n}^{\infty}\frac{1}{a_{n}}$.
     \par
     From Lemma \ref{lem:varbound}, $\mathbf{1}_{\llbracket{T_{n},\infty}\llbracket}\Delta{M^{n}}_{T_{n}}$ decomposes as $\mathbf{1}_{\llbracket{T_{n},\infty}\llbracket}\Delta{M^{n}}_{T_{n}}=D^{n}+C^{n}$ for some local martingale $D^{n}$ and some predictable process $C^{n}$ of integrable variation. Furthermore, Lemma \ref{lem:varbound} yields the bound,
     \begin{equation*}
         \int_{X}\mathrm{var}(C^{n})_{\infty}d\mu\leq\int_{X}\Vert{\Delta{M^{n}}_{T_{n}}}\Vert_{E}d\mu\lesssim 1.
     \end{equation*}
     By the Kadec-Pełczyński theorem, there exists a strictly positive sequence $\{b_{n}\}_{n}$, converging to $\infty$ in $\alpha[0,\infty)$, such that, after passing to a subsequence if necessary, $\{\mathrm{var}(C^{n})_{\infty}\wedge b_{n}\}_{n}$ is uniformly integrable.
     \par
     Define the predictable stopping times $\{V_{n}\}_{n}$ by
     \begin{equation*}
         V_{n}=\inf\{t:\mathrm{var}(C_{n})_{t}\geq b_{n}\}.
     \end{equation*}
     Thus, $M^{V_{k}-}=\mathbf{1}_{\llbracket{0,V_{k}}\llbracket}M+\mathbf{1}_{\llbracket{V_{k},\infty}\llbracket}M_{V_{k}-}$ is a martingale for any martingale $M$ and any $k\in\mathbb{N}$. It is easy to see from Markov's concentration inequality that $\mu(\{V_{n}<\infty\})\lesssim\frac{1}{b_{n}}$.
     \par
     Consider the martingale $N^{n}=\left((M^{n})^{T_{n}}-(\mathbf{1}_{\llbracket{T_{n},\infty}\llbracket}\Delta{M^{n}}_{T_{n}}-C^{n})\right)^{V_{n}-}$. Clearly,
     \begin{equation*}
         \left(N^{n}\right)^{\ast}_{\infty}\leq\left(M^{n}\right)^{\ast}_{\infty}\wedge a_{n}+\mathrm{var}(C^{n})_{\infty}\wedge b_{n},
     \end{equation*}
     implying that $\{N^{n}\}_{n}$ is relatively weakly compact in $H^{1}(\mu,E)$ (see Proposition \ref{prop:reflex-del-mey-yor}). Consider now the `singular part' $\{M^{n}-N^{n}\}_{n}$. An elementary calculation reveals,
     \begin{equation*}
         M^{n}-N^{n}=M^{n}-(M^{n})^{T_{n}\wedge V_{n}}+\Delta M^{n}_{T_{n}\wedge V_{n}}\mathbf{1}_{\llbracket{T_{n}\wedge V_{n},\infty}\llbracket}-(C^{n})^{V_{n-}}.
     \end{equation*}
     Let $L^{n}=M^{n}-(M^{n})^{T_{n}\wedge V_{n}}+\Delta M^{n}_{T_{n}\wedge V_{n}}\mathbf{1}_{\llbracket{T_{n}\wedge V_{n},\infty}\llbracket}$. Since 
     \begin{equation*}
         \lim_{n\to\infty}\mu(\{T_{n}\wedge V_{n}<\infty\})=0,
     \end{equation*}
     it follows that $\{L^{k}\}_{k}$ is a null sequence in the u.c.p. topology. The claim then follows from applying Lemma \ref{lem:pred-comp} to the sequence $\{R^{k}\}_{k}$ defined by $R^{k}=(M^{k}-N^{k})-L^{k}=-(C^{k})^{V_{k-}}$.
\end{proof}
Note that all variants of Komlós's theorem presented so far have assumed \textit{a priori} that the relevant sequences are bounded in $H^{1}(\mu,E)$. However, it is natural to inquire under what conditions does a Komlós-type result hold, without assuming $H^{1}(\mu,E)$ boundedness. We will now present such a result. Denote by $\mathscr{O}$ the optional $\sigma$-algebra on $[0,\infty)\times X$, which is generated by the right-continuous real-valued adapted processes.
\begin{theorem}\label{thm:unbound-komlos}
    Let $\{M^{n}\}_{n}$ be a sequence of $E$-valued càdlàg martingales, where $E$ is a reflexive Banach space. Suppose that
    \begin{equation*}
        \mathrm{co}\left\{\left(M^{n}\right)^{\ast}_{\infty}:n\in\mathbb{N}\right\},
    \end{equation*}
    is bounded in $L^{0}(\mu)$, and
    \begin{equation*}
        \left\{\Delta{M^{n}_{T}}:n\in\mathbb{N},T\textrm{ a stopping time}\right\},
    \end{equation*}
    is bounded in $L^{1}(\mu,E)$. Then there exists convex combinations $N^{i}\in\mathrm{co}\{M^{n}:n\geq i\}$ such that $\{N^{i}\}_{i}$ converges pointwise in $\sigma(E,E')$ to an $\mathscr{O}$-measurable process $N$ outside of an evanescent subset of $[0,\infty)\times X$.
\end{theorem}
\begin{remark}
We remark that the results of Theorem \ref{thm:unbound-komlos} are even relevant in the finite-dimensional case. Indeed, Theorem \ref{thm:unbound-komlos} refines (Theorem 2.6, \cite{martconvcomp}) by providing sufficient conditions under which one has convergence $\mu$-a.e., as opposed to merely in $L^{0}(\mu)$. This is especially important since convergence $\mu$-a.e. need not hold in the context of (Theorem 2.6, \cite{martconvcomp}), even under continuity assumptions (see Proposition 4.1, \cite{martconvcomp}). Furthermore, \cite{martconvcomp} only considers nonnegative martingales, while Theorem \ref{thm:unbound-komlos} deals with general martingales.
\end{remark}
The conditions in Theorem \ref{thm:unbound-komlos} are not very stringent. For context, one must understand that the functional-analytic (and convex-analytic) structure of $L^{0}(\mu)$ is quite poor. In particular, the convex hull of a bounded subset of $L^{0}(\mu)$ need not be bounded (see Example 1.2, \cite{kardzit} for a particularly ill-behaved example). As such, the condition is designed to eliminate some pathological examples.
\begin{proof}
    Since $K=\mathrm{co}\left\{\left(M^{n}\right)^{\ast}_{\infty}:n\in\mathbb{N}\right\}$ is convex, bounded in $L^{0}(\mu)$, and consists of nonnegative measurable functions, it follows from (Lemma 2.3, \cite{bipolar}) that there is a finite measure $\nu$ with $\nu\ll\mu\ll\nu$ such that
    \begin{equation*}
        \sup_{\zeta\in K}\int_{X}\vert{\zeta}\vert d\nu<\infty.
    \end{equation*}
    Thus, (Lemma 2.5, \cite{schruessdies}) and the triangle inequality yields convex combinations $N^{i}\in\mathrm{co}\{M^{n}:n\geq i\}$ and a nonnegative $\xi\in L^{0}(\mu)$ such that
    \begin{equation*}
        \sup_{i}\left(N^{i}\right)^{\ast}_{\infty}\leq\xi,
    \end{equation*}
    up to a $\nu$-null (equivalently, $\mu$-null) set. The singleton set $\{\xi\}\subset L^{0}(\mu)$ is bounded in $L^{0}(\mu)$, so that for each $\varepsilon>0$, there exists a $\delta_{\varepsilon}>0$ such that
    \begin{equation*}
        \mu(\{\xi\geq\delta_{\varepsilon}\})\leq\varepsilon.
    \end{equation*}
    Let $\{\varepsilon_{n}\}_{n}\subset(0,\infty)$ be a null sequence, and let $\{\delta_{\varepsilon_{n}}\}_{n}\subset(0,\infty)$ be as above; without loss of generality, we may assume that $\{\delta_{\varepsilon_{n}}\}_{n}$ is increasing. Define the sequence $\{T_{n}\}_{n}$ of stopping times by
    \begin{equation*}
        T_{n}=\inf\left\{t:\sup_{k}\left(N^{k}\right)^{\ast}_{t}\geq\delta_{\varepsilon_{n}}\right\}.
    \end{equation*}
    Clearly,
    \begin{equation*}
        \mu(\{T_{n}<\infty\})\leq\mu(\{\xi\geq\delta_{\varepsilon_{\varepsilon_{n}}}\})\leq\varepsilon_{n},
    \end{equation*}
    so that there is a measurable $X_{0}\in\mathscr{F}$ with $\mu(X\setminus X_{0})=0$ such that for each $\omega\in X_{0}$ there exists $n_{\omega}\in\mathbb{N}$ with $T_{m}(\omega)=\infty$ for all $m\geq n_{\omega}$.
    \par
    For each $n,i\in\mathbb{N}$,
    \begin{equation*}
        \left(N^{i}\right)^{\ast}_{T_{n}}\leq\delta_{\varepsilon_{n}}+\left\Vert{\Delta N^{i}_{T_{n}}}\right\Vert_{E},
    \end{equation*}
    by the triangle inequality. Thus, for each $n\in\mathbb{N}$,
    \begin{equation*}
        \sup_{i}\left\Vert{\left(N^{i}\right)^{T_{n}}}\right\Vert_{H^{1}(\mu,E)}<\infty.
    \end{equation*}
    By Theorem \ref{thm:ucp-komlos}, Theorem \ref{thm:del-mey-yor}, and a diagonalization procedure, we may pass to convex combinations (still denoted $\{N^{i}\}_{i}$) such that for all $n\in\mathbb{N}$,
    \begin{equation}\label{eq:kom-limit}
        \lim_{i\to\infty}(N^{i})^{T_{n}}=L^{n},
    \end{equation}
    in $\sigma(E,E')$ on $[0,\infty)\times X_{n}$ for some $X_{n}\in\mathscr{F}$ with $\mu(X\setminus X_{n})=0$, and some $\mathscr{O}$-measurable process $L^{n}$. Clearly, $L^{n}=L^{n+1}$ on $\llbracket{0,T_{n}}\rrbracket$. Define
    \begin{equation*}
        N=\mathbf{1}_{\llbracket{0,T_{1}}\llbracket}L^{1}+\sum_{n=2}^{\infty}\mathbf{1}_{\llbracket{T_{n-1},T_{n}}\llbracket}L^{n},
    \end{equation*}
    which converges absolutely up to an evanescent set. Let $\omega\in\bigcap_{n=0}^{\infty}X_{n}=\widetilde{X}$; then $T_{m}(\omega)=\infty$ for all $m\geq n_{\omega}$. Thus,
    \begin{equation*}
        N^{i}(\omega)=(N^{i})^{T_{n_{\omega}}(\omega)}(\omega),
    \end{equation*}
    where the right-hand side converges pointwise on $[0,\infty)$ in $\sigma(E,E')$ to $L^{n_{\omega}}(\omega)$ as $i\to\infty$ (see (\ref{eq:kom-limit})). Note that $N(\omega)=L^{n_{\omega}}(\omega)$ in this case. This proves the claim, as the complement of $[0,\infty)\times\widetilde{X}$ is evanescent, and $N$ is easily seen to be $\mathscr{O}$-measurable.
\end{proof}
\section*{Acknowledgments}
 The author would like to thank Professor V. Troitsky for organizing the University of Alberta seminar on functional analysis, where several productive conversations took place. The author would also like to thank Professor W. Schachermayer for carefully reading through the paper, and providing suggestions and advice. The author is additionally indebted to the two anonymous referees for their comments.
\printbibliography
\end{document}